\documentclass[reqno]{amsart}

\usepackage{amsmath}
\usepackage{amssymb}
\usepackage{amsfonts}

\usepackage{hyperref}
\usepackage{url}
\usepackage{mathrsfs}

\newtheorem*{thm}{Theorem}
\newtheorem{lem}{Lemma}

\begin{document}

\flushbottom

\title[On a Diophantine Inequality with Primes Yielding Square-Free Sums]{On a Diophantine Inequality with Primes Yielding Square-Free Sums with Given Numbers}

\author[T. Peneva]{Temenoujka P. Peneva}
\address{Faculty of Mathematics and Informatics, University of Plovdiv ``Paisii Hilendarski'', Plovdiv, Bulgaria}
\email{tpeneva@uni-plovdiv.bg}

\author[T. Todorova]{Tatiana L. Todorova}
\address{Faculty of Mathematics and Informatics, Sofia University ``St. Kliment Ohridski'', Sofia, Bulgaria}
\email{tlt@fmi.uni-sofia.bg}

\begin{abstract}
Let $\alpha\in \mathbb{R}\setminus\mathbb{Q}$ and $\beta\in \mathbb{R}$ be given. Suppose that $a_1,\ldots,a_s$ are distinct positive integers that do not contain a reduced residue system modulo $p^2$ for any prime $p$. We prove that there exist infinitely many primes $p$ such that the inequality $||\alpha p+\beta||<p^{-1/10}$ holds and all the numbers $p+a_1,\ldots,p+a_s$ are square-free.

\smallskip



\end{abstract}

\maketitle

\section{Introduction}
Let $r\ge 2$ be an integer. A natural number $n$ is called $r$-free if it is not divisible by the $r$th power of any prime $p$. In particular, 2-free numbers are also known as square-free numbers. 

Define $\mu_r(n)$ as the characteristic
function of the sequence of $r$-free numbers, i.e. $\mu_r(n)$ takes the value $1$ if $n$ is $r$-free, and $0$ otherwise. If $\mu$ denotes the M\"{o}bius function, it is easy to verify that 
\begin{equation*}
  \mu_r(n)=\sum\limits_{d^r|n}\mu (d).
\end{equation*}

Let $s\ge 2$ be an integer and $a_1,\ldots,a_s$ be distinct positive integers. The frequency of occurrence of systems of $r$-free numbers was first studied in 1936 by Pillai \cite{Pillai} for $r=2$, who established an asymptotic formula, with an error term $O(x/\log x)$, for the number of systems of square-free numbers $n+a_1,n+a_2,\ldots$, $n+a_{s}$ not exceeding $x$. This result was later generalized by Mirsky \cite{Mir}, \cite{Mir2} for any $r\ge 2$, who proved that for any $\varepsilon>0$,
\begin{equation}\label{Mirs}
  \sum\limits_{n\le x}\mu_r(n+a_1)\ldots\mu_r(n+a_s)=x\prod\limits_{p}\left(1
 -\frac{\nu(p^r)}{p^r}\right)
  +O\left(x^{\frac{2}{r+1}+\varepsilon}\right),
\end{equation}
where $\nu(p^r)$ is the number of distinct residue classes modulo $p^r$, represented by the numbers $a_1,\ldots,a_s$. 

Changa \cite{Changa} considered the case where $n$ is restricted to the set of prime numbers and obtained that for any $A>0$,
\begin{equation}\label{rezChanga}
  \sum\limits_{p\le x}\mu_{r}(p+a_1)\ldots\mu_{r}(p+a_s)=\pi(x)\prod\limits_{p}\left(1-\frac{\nu^*(p^r)}{\varphi(p^r)}\right)
  +O\left(\frac{x}{(\log x)^{A}}\right), 
\end{equation}
where $\varphi$ denotes the Euler function, and $\nu^*(p^r)$ is the number of distinct residue classes modulo $p^r$ that are co-prime with $p$, represented by the numbers $a_1,\ldots,a_s$. 

Observe that the infinite product in (1) (respectively, (2)) remains positive as long as, for any prime \( p \), the numbers \( a_1, \dots, a_s \) do not contain a complete (respectively, reduced) residue system modulo \( p^r \).

A more general problem was considered by Hablizel \cite{Habl}. 
For fixed $r_1,\ldots,$ $r_s \in \mathbb{N}$ satisfying
 $2\le r_1\le \ldots\le r_s$, he derived the asymptotic formula
\begin{equation*}
  \sum\limits_{p\le x}\mu_{r_1}(p+a_1)\ldots\mu_{r_s}(p+a_s)=\frac{x}{\log x} \prod\limits_{p}
  \left(1-\frac{D^*(p)}{\varphi(p^{r_s})}\right)
  +o\left(\frac{x}{\log x}\right),
\end{equation*}
where $D^*(p)$ is a computable function of the prime $p$, depending on the choice of the numbers $a_i$ and $r_i$.

Next, suppose that $\alpha$ is an irrational number and $\beta$ is any real number. A fundamental question in number theory concerns the validity of the Diophantine inequality 
\begin{equation}\label{1}
  ||\alpha p+\beta||<p^{-\theta}
\end{equation}
for infinitely many primes $p$, where, as usual, $||y||$ denotes the distance from $y$ to the nearest integer. 

In 1947, I.~M. Vinogradov \cite{Vin} first demonstrated that if $0<\theta<1/5$, then there exist infinitely many primes $p$ such that \eqref{1} holds. Subsequent research extended the range of the exponent $\theta$, with the most recent result, $0<\theta<1/3$, established by Matom\"{a}ki \cite{Mato2}. 

A natural variation of this problem involves restricting $p$ in inequality \eqref{1} to a specific subset of prime numbers (see, e.g., \cite{toltod}). In this paper, we take the set of primes $p$ for which $p+a_1,\ldots,p+a_s$ are square-free. 

\smallskip

We shall prove the following

\begin{thm} Let $\alpha$ be an irrational number and $\beta$ be a real number. Suppose $s\ge 2$ is an integer, and let $a_1<\ldots<a_s$ be positive integers that do not contain a reduced residue system modulo $p^2$ for any prime $p$. Then, for any $\theta<1/10$, there
exist infinitely many primes $p$ satisfying inequality \eqref{1} such that all the numbers $p+a_1,\ldots,p+a_s$ are
square-free.
\end{thm}

\noindent \textbf{Notation.} Let $x$ be a sufficiently large integer. Define
\begin{equation}\label{us2}
    \Delta =\Delta (x)=x^{-\theta },\quad \,K=\Delta ^{-1}\log ^2x,\quad \mbox{where }\theta<\frac{1}{10}.
\end{equation}
Throughout this paper, $p$ denotes a prime number. Instead of writing $m\equiv n \pmod k$
we use the shorthand notation $m\equiv n\,(k)$. For real $y$, we write $||y||$ for the distance from $y$ to the
nearest integer, $e(y)=\exp^{2\pi iy}$. As usual, $\mu(n)$, $\varphi(n)$, $\Lambda(n)$, and $\tau_k(n)$ denote the M\"{o}bius function, the Euler function, the von Mangoldt function, and the $k$th divisor function, respectively; $\tau(n)=\tau_2(n)$. The function $\nu_p(n)$ is defined such that $\nu_p(n)=k$ if $p^k|n$ but $p^{k+1}\nmid n$.

If $X$ and $Y$ are positive numbers, the notation \( X \asymp Y \) signifies that \( X \ll Y \ll X \). Furthermore, \( n \sim X \) indicates that \( n \) runs through some subinterval of \( (X, 2X] \), although the precise endpoints may vary depending on the context. The symbol \( \varepsilon \) denotes an arbitrarily small positive number, whose value may change from one occurrence to another. This convention allows for the use of inequalities such as \( X^\varepsilon \log X \ll X^\varepsilon \).

\section{Auxiliary results}

Before launching the proof of Theorem 1, we prepare the ground with some auxiliary results for the reader's convenience.

The first two statements correspond to Lemmas 8 and 9 of Mennema [8]. They provide average bounds for the divisor function over square-free numbers.

\begin{lem}\label{Lem1Mennema}There exists $C_1>1$ such that for all integer $k\ge 2$ 
and for all real $x\ge 1$,
\begin{equation*}
  \sum\limits_{n\le x}\mu ^2(n)\tau_k(n)\le C_1^kx(\log x)^{k-1}.
\end{equation*}
\end{lem}

\begin{lem}\label{Lem2Mennema} There exists $C_2>1$ such that for all integer $k\ge 2$
and for all real $x\ge 1$,
\begin{equation*}
  \sum\limits_{d> x}\frac{\mu ^2(d)\tau_k(d)}{d^2}\le \frac{C_2^k(2k-2+\log x)^{k-1}}{x}.
\end{equation*}
\end{lem}

Let $n,\varpi\in \mathbb{N}$. Following the notation of Mennema \cite[\S3]{Mennema}, we write 
\begin{equation}\label{mutildew}
  \mu (n)=\mu_{\varpi}(n)\widetilde{\mu}(n),
\end{equation}
where
\begin{equation}\label{defmutilde}
  \mu_{\varpi}(n)=\mu\Bigg(\prod\limits_{p|\varpi}p^{\nu_p(n)}\Bigg),\qquad
  \widetilde{\mu}(n)=\mu\Bigg( \prod\limits_{p\nmid \varpi}p^{\nu_p(n)} \Bigg). 
\end{equation}

The following two lemmas are Lemma 3.3 and Lemma 3.4 from Mennema \cite{Mennema}. 

\begin{lem} \label{Mennema1} Let $n, m, \varpi \in \mathbb{N}$ be such that $n\equiv m \pmod {\varpi^2}$, and let the function $\mu_{\varpi}$ be as in \eqref{defmutilde}. Then
$\mu_{\varpi}(n) = \mu_{\varpi} (m)$. 
\end{lem}

\begin{lem} \label{Mennema2} Let $n, \varpi \in \mathbb{N}$, and let $\widetilde{\mu}$ be as in \eqref{defmutilde}. Then
\begin{equation*}
  \widetilde{\mu}^2(n)=\sum\limits_{\substack{d^2|n \\ (d,\,\varpi)=1}}\mu (d).
\end{equation*}
\end{lem}

From this point onward, we put $\varpi$ to be the constant
\begin{equation}\label{lcmw}
  \varpi=\prod\limits_{p\le (a_s-a_1)^{1/2} }p.
\end{equation}

The following lemma is Lemma 3.5 from \cite{Mennema}.

\begin{lem} \label{new}
Let $n,a_1,\ldots,a_s$ be positive integers, and let $a_1<\ldots<a_s$. If $d_i^2|n+a_i$, $d_j^2|n+a_j$, and $(d_id_j,\,\varpi)=1$, then $(d_i,\,d_j)=1$ for all $i\ne j$. 
\end{lem}

The proof of our Theorem will depend on estimates of exponential sums. The following statement is a direct consequence of Lemma 4 in \cite[Chapter 6, \S2]{Kar}.

\begin{lem} \label{Trsum1}
Let $ X\ge 1 $ and $\alpha$ be real numbers, $a$, $d\in \mathbb{Z}$, $d\ge 1$. Then
\begin{equation*}
   \left|
\mathop{\sum}_{\substack{
                      n\le X \\
                      n\equiv a \,(d) }}
e(\alpha n)\right|
           \ll \min \left\{\frac{X}{d},
                 \frac{1} { ||\alpha d||} \right\}.
\end{equation*}
\end{lem}

Furthermore, suppose that $\alpha$ is a real number with a rational approximation $a/q$ satisfying
\begin{equation}\label{alfa}
  \left|\alpha- \frac{a}{q}\right|<\frac{1}{q^2},\quad \mbox{where} \;  (a, q) = 1 \;\mbox{and}\; q\ge 1.
\end{equation}

\smallskip

The following lemma is a well-known estimate of Vaughan \cite[Chapter 2, \S2.1]{Va}. 

\begin{lem} \label{Trsum2}
Suppose that $ X,\,Y\ge 1$ are real numbers, and that $\alpha$ is a real number satisfying \eqref{alfa}. Then
\begin{equation*}
  \sum_{n\le X} \min \left\{ \frac{XY}{n}, \frac{1} { ||\alpha n|| } \right\}
    \ll  XY
      \left( \frac{1} {q} +
                \frac{1} {Y} +
                \frac { q }{ XY } \right)\log (2Xq).
\end{equation*}
\end{lem}

The next lemma is a consequence of Matom\"{a}ki's result \cite[Lemma 8]{Mato1}. 

\begin{lem}\label{Mat} 

Suppose that $x,\,M,\,J\in \mathbb{R}^{+}$, $\mu,\,\zeta\in\mathbb{N}$, and that $\alpha$ is a real number satisfying \eqref{alfa}. Then for any $\varepsilon>0$,
\begin{multline*}
  \sum\limits_{m\sim M}\tau_{\mu}(m)\sum\limits_{j\sim J}\tau_{\zeta}(j)
\min\left\{\frac{x}{m^2j},\,\frac{1}{||\alpha m^2j||}\right\}\\
\ll x^{\varepsilon}\left(MJ+\frac{x}{M^{3/2}}+\frac{x}{Mq^{1/2}}+
\frac{x^{1/2}q^{1/2}}{M}\right).
\end{multline*}

\end{lem}

The following statement is \cite[Lemma 8]{TLT}. 

\begin{lem}\label{TT} Suppose that $x,\,M,\,J\in \mathbb{R}^{+}$, $\mu,\,\zeta\in \mathbb{N}$,
and that $\alpha$ is a real number satisfying
\eqref{alfa}. Then for any $\varepsilon>0$,  
\begin{multline*}
  \sum\limits_{m\sim M}\tau_{\mu}(m)\sum\limits_{j\sim J}\tau_{\zeta}(j)\min\left\{\frac{x}{m^4j},\,\frac{1}{||\alpha m^4j||}\right\}
\\ \ll x^{\varepsilon}\left(MJ+\frac{x}{M^{25/8}}+\frac{x}{M^3q^{1/8}}+
\frac{x^{7/8}q^{1/8}}{M^3}\right).
\end{multline*}
\end{lem}

\section{Proof of the Theorem}

We start by observing that there exists a periodic function $\chi$ with period 1 such that
\begin{equation}\label{hi}\nonumber
\begin{gathered}
    0<\chi (t) \le 1 \; \mbox { for }\;  -\Delta< t< \Delta,\\
    \chi (t) =0 \; \mbox { for }\;  \Delta \le t\le 1-\Delta,
\end{gathered}
\end{equation}
and $\chi(t)$ admits a Fourier expansion
\begin{equation}\label{hi1}
    \chi (t)= \Delta +\Delta\sum\limits_{|k|>0}g(k)e(kt),
\end{equation}
where the Fourier coefficients satisfy
\begin{equation}\label{hi2}
  g(k)\ll 1  \;\mbox{ for all }\; k\ne 0,\qquad
    \Delta\sum\limits_{|k|>K}|g(k)|\ll x^{-1}.
\end{equation}
The existence of such a function is a consequence of a lemma of Vinogradov (see \cite[Chapter 1, \S2]{Kar}).

Consider the sum
\begin{equation*}
   \Gamma (x)= \sum\limits_{p\sim x }\chi (\alpha p+\beta)\mu ^2(p+a_1)\ldots\mu ^2(p+a_s).
\end{equation*}
To prove our theorem, it suffices to determine the
constant $\theta$ such that
there exists a sequence of positive integers $\{x_j\}_{j=1}^\infty$ satisfying 
\begin{equation}\label{end}
    \lim\limits_{j\to\infty}x_j=\infty \qquad
\end{equation}
and
\begin{equation}\label{end2}
  \Gamma(x_j)\ge  \frac{C\Delta (x_j)\,x_j}{\log x_j},\,\quad j=1,2,3,\ldots
\end{equation}
with some absolute constant $C>0$. 
    
Applying the Fourier expansion \eqref{hi1} and the inequalities \eqref{hi2}, we obtain
\begin{equation}\label{RavGama1}
    \Gamma (x)=\Delta\left(\Gamma_1(x)+\Gamma_2(x)\right)+O(1),
\end{equation}
where
\begin{equation*}
 \Gamma_1(x) =\sum\limits_{p\sim x }\mu ^2(p+a_1)\ldots\mu ^2(p+a_s),
\end{equation*}
and
\begin{equation} \label{Gamma_2}
\Gamma_2(x) =
   \sum\limits_{0<|k|\le K}c(k)
   \sum\limits_{p\sim x}\,
   \mu ^2(p+a_1)\ldots\mu ^2(p+a_s) e(\alpha k p),  \,
\end{equation}
with $c(k):=g(k)e(\beta k)$, satisfying
\begin{equation}\label{c(k)}
    c(k)\ll 1  \;\mbox{ for all }\; k\ne 0. 
\end{equation}

Consider the sum $\Gamma_1(x)$. Changa's asymptotic formula (\ref{rezChanga}) yields that for any $A>0$,
$$
  \Gamma_1(x)=\mathfrak{S}\left( \pi (2x)-\pi(x)\right)
  +O\left( \frac{x}{(\log x)^{A}} \right),
$$
where 
$$
  \mathfrak{S}= \prod_{p} \left( 1-\frac{\nu^*(p^2)}{p(p-1)} \right)
$$ 
is the infinite product in \eqref{rezChanga} for $r=2$. Observe that $\mathfrak{S}>0$, since for the given $a_1,\ldots,a_s$, every factor in $\mathfrak{S}$ is positive, and the factor corresponding to $p$ is at least $1-s/(p(p-1))$ for all sufficiently large values of $p$. 
According to Rosser and Schoenfeld's classic estimate~\cite[Corollary 3]{Rosser}, for  $x\ge 20.5$,
$$
  \pi(2x)-\pi(x)>\frac{3x}{5\log x}.
$$
Thus, it follows that
\begin{equation}\label{Gamma1}
  \Gamma_1(x)> \frac{\mathfrak{S}x}{2\log x},
\end{equation}
for sufficiently large $x$. 

The estimate of $\Gamma_2(x)$ is postponed until Section 4.

\vskip 1pt

Now, let $\{q_j\}_{j=1}^\infty$ be a sequence
of values of $q$ that satisfy (\ref{alfa}). In view of \eqref{us2} and the estimate of $\Gamma_2(x)$ obtained in Section \ref{conclusion}, we define a sequence $\{x_j\}_{j=1}^\infty$  by setting 
$$
x_j=q_j^{20/7}, \quad j=1,2,\ldots\,.
$$
Condition \eqref{end} is clearly satisfied. Furthermore, for a sufficiently small $\varepsilon>0$ and any $A>0$, we have
\begin{equation}\label{GamaOkonch}
  \Gamma_2(x_j)\ll x_j^{9/10+\varepsilon}K\ll \frac{x_j
  }{(\log x_j)^{A}},\quad i=1,2,\ldots\,.
\end{equation}
Using (\ref{RavGama1}), (\ref{Gamma1}) and (\ref{GamaOkonch}), we deduce the estimate \eqref{end2}  
with some absolute constant $C<\mathfrak{S}/2 $, thus completing the proof of the Theorem.  

\section{The estimation of the sum $\Gamma_2$} \label{g2} 

In this section, we estimate the sum $\Gamma_2(x)$, as defined in $\eqref{Gamma_2}$.

\subsection{Preparation} We begin by adapting Mennema's argument from \cite[\S3]{Mennema}. For $n\in\mathbb{N}$, define 
\begin{equation}\label{fmdef}
  f(n)=\mu_{\varpi}^2(n+a_1)\ldots\mu_{\varpi}^2(n+a_s), 
\end{equation}
where $\mu_{\varpi}$ and $\varpi$ are given in \eqref{defmutilde} and \eqref{lcmw}, respectively. 
Obviously, $f(n)=1$ if $(n+a_i,\varpi^2)$ is square-free for all $i$, and $f(n)=0$ otherwise. In all cases, 
\begin{equation}\label{f(t)}
    0\le f(n)\le 1 \quad \mbox{for all} \quad n\in \mathbb{N}.
\end{equation}

Now, using \eqref{mutildew},  \eqref{fmdef}, and Lemma \ref{Mennema2}, we can express
\begin{align*}
    \Gamma_2(x) & =
   \sum\limits_{0<|k|\le K}c(k)
   \sum\limits_{p\sim x}\,
   f(p)\,\tilde{\mu}^2(p+a_1)\ldots\tilde{\mu}^2(p+a_s) e(\alpha k p)\\
   &=\sum\limits_{0<|k|\le K}c(k)  
   \sum\limits_{\substack{p\sim x}} f(p)e(\alpha k p )\sum\limits_{\substack{d_i^2| p+a_i\\(d_i,\varpi)=1\\i=1,\ldots,s}} \mu(d_1)\ldots \mu(d_s).
\end{align*}
Next, we decompose the sum over $p$ into sums over residue classes modulo $\varpi^2$.  Noting that, by Lemma \ref{Mennema1}, $p\equiv t\,(\varpi^2)$ implies $f(p)=f(t)$, and changing the order of summation, we obtain 
\begin{align}
   \Gamma_2(x)& = \sum\limits_{1\le t\le \varpi^2 }f(t)\!\sum\limits_{0<|k|\le K}c(k) \sum\limits_{\substack{p\sim x\\p\equiv t\,(\varpi^2)}} e(\alpha kp)\sum\limits_{\substack{d_i^2|p+a_i\\(d_i,\varpi)=1\\i=1,\ldots,s}} \mu(d_1)\ldots \mu(d_s)  \notag\\
   & = 
   \sum\limits_{1\le t\le \varpi^2 }f(t) \sum\limits_{0<|k|\le K}c(k) \sum\limits_{\substack{1\le d_i\le X_i\\ (d_i,\varpi)=1\\ i=1,\ldots,s}}\! \mu (d_1)\ldots\mu(d_s) 
   \! \sum\limits_{\substack{p\sim x\\ p\equiv t\,(\varpi^2)\\p\equiv -a_i\,(d_i^2)\\i=1,\ldots,s}}e(\alpha kp),\notag
\end{align}
with
\begin{equation}\label{X_i}\nonumber
   X_i:=(2x+a_i)^{1/2}, \qquad X_i\asymp x^{1/2},\qquad i=1,\ldots,s.
\end{equation}
Note that if there exists an index \( i \) such that \( d_i \) and \( a_i \) share a common prime factor \( p' \), then the condition \( p \equiv -a_i \,(d_i^2) \) implies that $p'$ divides $p$, which further leads to \( p'=p \) and \( p' \sim x \). However, this contradicts \( p'\le d_i \ll x^{1/2} \). Consequently, we may assume $ (d_i,a_i)=1$ , and hence $(d_i,\varpi a_i)=1$ for $i=1,\ldots,s$. Following the same line of reasoning, we conclude $(t,\varpi)=1$. By using \eqref{f(t)}, we obtain
\begin{align}
   \Gamma_2(x)& \le \varpi^2\max\limits_{1\le t\le \varpi^2\atop{(t,\varpi)=1}}\left| \sum\limits_{0<|k|\le K}c(k) \sum\limits_{\substack{1\le d_i\le X_i\\ (d_i,\varpi a_i)=1\\ i=1,\ldots,s}}\! \mu (d_1)\ldots\mu(d_s) 
   \! \sum\limits_{\substack{p\sim x\\ p\equiv t\,(\varpi^2)\\p\equiv -a_i\,(d_i^2)\\i=1,\ldots,s}}e(\alpha kp)\right|.\notag
\end{align}

By partial summation, we have
\begin{equation}\label{Gamma_2-2}
    \Gamma_2(x) \ll (\log x)^{-1} \max\limits_{\substack{x<y\le 2x\\ 1\le t\le \varpi^2\\(t,\varpi)=1}}\big|\Gamma_3(y) \big|
    + O(K x^{1/2}\log x), 
\end{equation}
where 
\begin{equation}\nonumber
   \Gamma_3(y)=\sum\limits_{0<|k|\le K}c(k) \sum\limits_{\substack{1\le d_i\le X_i\\ (d_i,\varpi a_i)=1\\ i=1,\ldots,s}} \mu (d_1)\ldots\mu(d_s) 
    \sum\limits_{\substack{n\sim y\\ n\equiv t\,(\varpi^2)\\n\equiv -a_i\,(d_i^2)\\i=1,\ldots,s}} \Lambda(n)e(\alpha k n).
\end{equation}

\noindent It is evident that the conditions imposed on the summation over $d_i$'s imply $d_1\ldots d_s$ $\in[1,X_1\ldots X_s]$. Thus, we decompose this interval as: 
$$ 
  [1,X_1\ldots X_s]=\mathcal{I}_1\cup \mathcal{I}_2\cup \mathcal{I}_3,
$$
where the intervals are defined as:
$$
  \mathcal{I}_1 = (y^{1/2},X_1\ldots X_s],\qquad \mathcal{I}_2 =(y^{1/5},y^{1/2}], \qquad \mathcal{I}_3 = [1,y^{1/5}].
$$
Consequently, we express the sum $\Gamma_3(y)$ as the sum of three components:
\begin{equation}\label{new6}
\Gamma_3 (y) = \mathcal{U}_1+\mathcal{U}_2+\mathcal{U}_3,
\end{equation}
with each $\mathcal{U}_i$ defined by
\begin{equation}\label{V_d}
   \mathcal{U}_i= \sum\limits_{0<|k|\le K}c(k) \sum\limits_{\substack{1\le d_i\le X_i\\ (d_i,\varpi a_i)=1\\ i=1,\ldots,s\\d_1\ldots d_s\in \mathcal{I}_i}}\! \mu (d_1)\ldots\mu(d_s) 
   \! \sum\limits_{\substack{n\sim y\\ n\equiv t\,(\varpi^2)\\n\equiv -a_i\,(d_i^2)\\i=1,\ldots,s}} \!\Lambda(n)e(\alpha k n).
\end{equation}

The remainder of this paper is devoted to estimating these quantities. 

\subsection{Estimating \boldmath $\mathcal{U}_1$\unboldmath}  To estimate \/$\mathcal{U}_1$, we follow Mennema's approach from \cite[\S3]{Mennema}.

First, we impose the ordering constraint $d_1\le \ldots \le d_s$, which introduces an additional factor of $s!$ in our estimate. For integers $r$ satisfying $2\le r\le s$, we define
\begin{equation*}
  \alpha_r=\frac{s-r+2}{s-r+1},\qquad A_r=y^{\frac{1}{1+\alpha_r}}=y^{\frac{s-r+1}{2s-2r+3}},
\end{equation*}
and the corresponding set
\begin{equation*}
\begin{split}
     &\mathcal{D}_r=\bigg\{(d_1,\ldots,d_s)\in \mathbb{N}^s\;:
\end{split}
\quad\;
  \begin{split}
     & d_1\le \ldots\le d_s, \\
     & d_r\ldots d_s\le A_r, 
  \end{split}
\quad\;
  \begin{split}
     & (d_i,\varpi)=1 \mbox{ for all } i \; \\
     & d_{r-1}\ldots d_s>A_{r-1}
  \end{split}
\bigg\}.
\end{equation*}
Notice that, according to Lemma \ref{new}, we have $(d_i,d_j)=1$ whenever $i\ne j$. Since $d_1\ldots d_s> y^{1/2}$, there exists an $r$ such that $(d_1,\ldots,d_s)\in \mathcal{D}_r$ (see \cite[Remark 3.7]{Mennema}). 
On writing 
$$
  \delta=d_{r-1}\ldots d_s,
$$
we deduce $A_{r-1}<\delta\le A_r^{\alpha_r}$ and $(\delta,\varpi )=1$ (see \cite[p.\,20]{Mennema}). Therefore, the contribution of $(d_1,\ldots,d_s)\in \mathcal{D}_r$ to the sum \/ $\mathcal{U}_1$ is
\begin{align*}
   &\ll y^\varepsilon\sum\limits_{0<k\le K}
   \sum\limits_{\substack{A_{r-1}<\delta\le A_r^{\alpha_r} \\ (\delta,\,\varpi)=1}} \mu ^2(\delta)\tau_{s-r+2}(\delta)
   \sum\limits_{\substack{n\sim y\\ n\equiv -a_i\,(d_i^2)\\i=r-1,\ldots,s}}
   \sum\limits_{\substack{1\le d_i\le X_i\\ (d_i,\varpi a_i)=1\\ d_i^2|n+a_i\\ i=1,\ldots,r-2}}1\\
   &\ll y^{\varepsilon}K
   \sum\limits_{\substack{A_{r-1}<\delta\le A_r^{\alpha_r}\\(\delta,\,\varpi)=1}}\mu ^2(\delta)\tau_{s-r+2}(\delta)
   \max\limits_{ 1\le a\le \delta^2}\sum\limits_{\substack{n\sim y\\n\equiv a\,(\delta^2)} }
   \tau(n+a_{1})\ldots\tau(n+a_{r-2})\\
   &\ll y^{\varepsilon}K\sum\limits_{\substack{A_{r-1}<\delta\le A_r^{\alpha_r}\\(\delta,\,\varpi)=1} }
   \mu ^2(\delta)\tau_{s-r+2}(\delta)\left(\frac{y}{\delta^2}+1\right),
\end{align*}
where we have used \eqref{c(k)} and the well-known estimate
\begin{equation}\label{Tau}
\tau_k(n)\ll_{k, \varepsilon} n^\varepsilon.
\end{equation}
Applying Lemmas \ref{Lem1Mennema} and \ref{Lem2Mennema}, we deduce
\begin{equation*}
  \mathcal{U}_1\ll y^{\varepsilon}K \sum_{r=2}^s  \left( A_r^{\alpha_r}(\log A_r^{\alpha_r})^{s}+\frac{y(2s+\log A_{r-1})^{s}}{A_{r-1}}\right).
\end{equation*}
For $2\le r\le s$, it is easy to verify that
$$
  A_r^{\alpha_r}=y^{\frac{s-r+2}{2s-2r+3}}\le y^{2/3}, \qquad 
  \frac{y}{A_{r-1}}< \frac{y}{A_{r}}=y^{\frac{s-r+2}{2s-2r+3}}\le y^{2/3}.
$$
Thus, we conclude that
\begin{equation}\label{veryLargeA0}
  \mathcal{U}_1\ll y^{2/3+\varepsilon}K.
\end{equation}

\subsection{Estimating \boldmath $\mathcal{U}_2$ \unboldmath} We estimate the sum $\mathcal{U}_2$, as defined in \eqref{V_d}, fairly straightforwardly. Setting $\tilde{d}=d\varpi$ and applying \eqref{c(k)} and \eqref{Tau}, we obtain
\begin{align}\label{BigA0}
  \mathcal{U}_2&\ll y^\varepsilon \sum\limits_{0<k\le K}\;\;
  \sum\limits_{\varpi y^{1/5}< \tilde{d}\le \varpi y^{1/2}}\tau_{s+1}(\tilde{d})\max\limits_{ 1\le a\le \tilde{d}^2}   
  \sum\limits_{\substack{n\sim y\\ n\equiv a\,(\tilde{d}^2)} } 1
\notag \\
&\ll y^{\varepsilon} K
  \sum\limits_{\varpi y^{1/5}\le \tilde{d}\le \varpi y^{1/2}}
  \left(\frac{y}{{\tilde{d}}^2}+1\right)\notag \\
&\ll y^{4/5+\varepsilon}K.
\end{align}

\subsection{Estimating \boldmath $\mathcal{U}_3$ \unboldmath}

By a dyadic decomposition of the summation ranges, we write $\mathcal{U}_3$, defined in \eqref{V_d}, as a sum of $O(\log ^{s+1}x)$ sums of the type
\begin{equation*}
   W=\sum\limits_{k\sim K_0}c(k)
   \sum\limits_{\substack{d_i\sim D_i\\ (d_i,\varpi a_i)=1\\ i=1,\ldots,s\\d_1\ldots d_s \in \mathcal{I}_3}}\mu (d_1)\ldots\mu (d_s)
   \sum\limits_{\substack{n\sim y\\ n\equiv t\,(\varpi^2)\\ n\equiv -a_i\,(d_i^2)\\i=1,\ldots,s}} \Lambda (n)e(\alpha k n),
\end{equation*}
where
\begin{equation}\label{Ogranichenia}
   1\le K_0\ll K,\qquad   1\le D_i\ll X_i, \qquad i=1,\ldots,s.
\end{equation}
Define 
$$
D=D_1\ldots D_s.
$$
The conditions \( d_i \sim D_i \) for \( i = 1, \dots, s \) imply that \( D \leq d_1 \dots d_s \leq 2^s D \). Since \( d_1 \dots d_s \in \mathcal{I}_3\), we deduce
\begin{equation}\label{D}
  1\le D \le y^{1/5}.
\end{equation}
With this inequality in mind, we proceed to replace the condition $d_1\ldots d_s\in \mathcal{I}_3$  with $d_1\ldots d_s\asymp D$ in the sum over $d_i$'s.

Applying Heath-Brown's identity \cite{Heath} with parameters
\begin{equation*}\label{uvw}
u=2^{-7}y^{\frac{1}{5}},\quad v=2^7y^{\frac{1}{3}},\quad w=y^{\frac{2}{5}},
\end{equation*}
we decompose the sum $W$ as a linear combination of $O(\log ^6x)$
sums of type I and type II. The type I sums are
\begin{equation*}
    W_1=\sum\limits_{k\sim K_0}c(k)
   \sum\limits_{\substack{d_i\sim D_i\\ (d_i,\varpi a_i)=1\\i=1,\ldots,s\\ d_1\ldots d_s\asymp D }} \mu (d_1)\ldots\mu (d_s)
   \sum_{m\sim M} a(m)\sum\limits_{\substack{ \ell\sim L \\ m\ell \equiv t (\varpi^2)\\m\ell \equiv -a_i (d_i^2)\\i=1,\ldots,s}}
   \! a(m)e(\alpha m\ell k),
\end{equation*}
and
\begin{equation*}
    W'_1=\sum\limits_{k\sim K_0}c(k)
   \sum\limits_{\substack{d_i\sim D_i\\ (d_i,\varpi a_i)=1\\i=1,\ldots,s\\d_1\ldots d_s\asymp D}} \mu (d_1)\ldots\mu (d_s)
   \sum_{m\sim M} a(m) \sum\limits_{\substack{ \ell\sim L \\ m\ell \equiv t (\varpi^2)\\m\ell \equiv -a_i (d_i^2)\\i=1,\ldots,s}}
   e(\alpha m\ell k) \log \ell ,
\end{equation*}
where
\begin{equation}\label{Usl1}
 ML\asymp y, \quad  L\ge w, \quad  a(m) \ll  y^{\varepsilon}.
\end{equation}
The type II sums are
\begin{equation*}
    W_2= \sum\limits_{k\sim K_0}c(k)
   \sum\limits_{\substack{d_i\sim D_i\\ (d_i,\varpi a_i)=1\\i=1,\ldots,s\\ d_1\ldots d_s\asymp D}}\mu (d_1)\ldots\mu (d_s)
    \sum_{\ell\sim L}b(\ell)\sum\limits_{\substack{ m\sim M \\ m\ell \equiv t (\varpi^2)\\m\ell \equiv -a_i (d_i^2)\\i=1,\ldots,s}}
   a(m) b(\ell) e(\alpha m\ell k),
\end{equation*}
where
\begin{equation}\label{Usl3}
ML\asymp y, \quad u\le L\le v, \quad a(m), \; b(\ell) \ll  y^{\varepsilon}.
\end{equation}

Recalling \eqref{D}, in Sections \ref{new8} and \ref{new9} we establish bounds for the sums $W_1$, $W_1'$ and $W_2$, and hence for $W$, under the 
assumptions $y^{3/20}\le D\le y^{1/5}$ 
and $1\le D\le y^{3/20}$, respectively.

\subsubsection{Estimate of \/ $W$ in the case $y^{3/20} \le D\le y^{1/5}$}\label{new8}
Let us denote $\tilde{d}=\varpi \,d_1\ldots d_s $. Since $(t,\varpi)=1$, $(d_i,\,\varpi a_i)=1$, and $(d_i,d_j)=1$ for $i\ne j$, the Chinese Remainder Theorem implies that the system of congruences $n\equiv t\,(\varpi^2)$, $n\equiv a_i\,(d_i^2),\,i=1,...,s$, has a unique solution $a$ modulo $\tilde{d}^2$, with $(a,\tilde{d})=1$. Thus, using estimates  \eqref{c(k)}, \eqref{Tau}, and \eqref{Usl1}, we obtain
\begin{equation} \label{W1mkd}
  {W}_1 \ll y^{\varepsilon}\sum\limits_{k\sim K_0}
   \;\sum\limits_{\tilde{d}\sim \varpi D }\tau_{s+1}(\tilde{d})
    \;\sum\limits_{m\sim M}\;\max\limits_{\substack{1\le a\le \tilde{d}^2\\(a,\tilde{d})=1}}
    \left|\sum\limits_{\substack{ \ell\sim L \\ m\ell \equiv a (\tilde{d}^2)}}
e(\alpha m\ell k)\right|\notag.
\end{equation}
It is clear that if $(m,\tilde{d})>1$, then for any integer $a$ satisfying $(a,\tilde{d})=1$, the congruence $m\ell\equiv a\,(\tilde{d}^2)$ admits no solution for $\ell$. Conversely, if $(m,\tilde{d})=1$, the congruence $m\ell\equiv a\,(\tilde{d}^2)$ has a unique solution for $\ell$ modulo $\tilde{d}^2$. Applying  Lemma~\ref{Trsum1}, we obtain
\begin{equation}\label{W_1}
    W_1 \ll y^{\varepsilon}\sum\limits_{k\sim K_0}\;\sum\limits_{\tilde{d}\;\sim \varpi D}
  \;\sum\limits_{m\sim M}\min\bigg\{\frac{yK_0}{mk \tilde{d}^2},\,\frac{1}{||\alpha mk\tilde{d}^2||}\bigg\}.
\end{equation}
By setting $h=mk$ and employing (\ref{Ogranichenia}) and \eqref{Usl1}, along with Lemma \ref{Mat}, we derive
\begin{align}
  W_1&\ll y^{\varepsilon}\sum\limits_{h\sim MK_0}\sum\limits_{\tilde{d}\sim \varpi D}
  \min\bigg\{\frac{yK_0}{h\tilde{d}^2},\,\frac{1}{||\alpha h\tilde{d}^2||}\bigg\}\notag \\
  &\ll y^{\varepsilon}\left(y^{4/5}K+\frac{y^{17/20}K}{q^{1/2}}+y^{7/20}K^{1/2}q^{1/2}\right).
  \label{W1intermediate}
\end{align}


We get the same estimate for $W_1'$ by partial summation. 

To estimate $W_2$,  we apply the Cauchy-Schwarz inequality repeatedly, and appeal to \eqref{c(k)}, \eqref{Tau}, and \eqref{Usl3} to deduce
\begin{equation}\label{W2malkoAsquare}
  W_2^2\ll y^{\varepsilon}\left(y^{9/5}K^2 +W_{21}\right),
\end{equation}
where
\begin{align*}
  W_{21} & \ll MDK\sum\limits_{k\sim K_0}
  \sum\limits_{\substack{d_i\sim D_i\\ (d_i,\varpi a_i)=1\\ i=1,\ldots,s\\d_1\ldots d_s\asymp D}}
  \sum\limits_{\substack{\ell_1,\ell_2\sim L\\\ell_1 \ne \ell_2}}  b(\ell_1)b(\ell_2)  
  \sum\limits_{\substack{m\sim M \\ m\ell_j\equiv t (\varpi^2) \\ m\ell_j\equiv -a_i (d_i^2)\\ j=1,2 \\i=1,\ldots,s}} e(\alpha m(\ell_1-\ell_2)k)\\
   & \ll y^{\varepsilon}MDK
  \sum\limits_{k\sim K_0}\sum\limits_{\tilde{d}\sim \varpi D}\tau_{s+1}(\tilde{d})
  \sum\limits_{\substack{\ell_1,\ell_2\sim L\\\ell_1 \ne \ell_2}} 
  \max\limits_{\substack{1\le a\le \tilde{d}^2\\ (a,\tilde{d})=1}}\left|
  \sum\limits_{\substack{m\sim M\\ m\ell_j \equiv a (\tilde{d}^2)\\ j=1,2}}e(\alpha m(\ell_1-\ell_2)k)\right|,
\end{align*}
by reasoning analogously to the estimation of the sum $W_1$. From $\ell_1\not= \ell_2$ and $m\ell_j\equiv a(\tilde{d}^2)$ for $j=1,2$, we have $\ell_1=\ell_2+z\tilde{d}^2$.
By setting $h=zk$, we have
\begin{align}\label{W21square}
  W_{21}&\ll y^{\varepsilon}MDK\sum\limits_{k\sim K_0}\sum\limits_{\tilde{d}\sim \varpi D}
  \sum\limits_{\ell_2\sim L}\sum\limits_{z\le L/\tilde{d}^2} \max\limits_{\substack{1\le a\le \tilde{d}^2\\ (a,\tilde{d})=1}}\left|
  \sum\limits_{\substack{m\sim M\\ m\ell_2\equiv a(\tilde{d}^2)}} e(\alpha mt\tilde{d}^2k)
  \right|\notag.
\end{align}
As discussed previously, the congruence $m\ell_2\equiv a\,(\tilde{d}^2)$ either admits no solutions for $m$ or exactly one solution for $m$ modulo $\tilde{d}^2$. Consequently, by using Lemma \ref{Trsum1}, we obtain
\begin{align}
  W_{21}&\ll y^{\varepsilon}MLDK\sum\limits_{k\sim K_0}\sum\limits_{\tilde{d}\sim \varpi D}
  \sum\limits_{z\le L/\tilde{d}^2}
  \min\left\{
  \frac{yK_0}{\tilde{d}^4tk},\,\frac{1}{||\alpha \tilde{d}^4tk||}
  \right\}\notag\\
  &\ll y^{1+\varepsilon}DK\sum\limits_{h\ll \frac{LK_0}{\varpi^2D^2}}
  \sum\limits_{\tilde{d}\sim \varpi D}
  \min\left\{\frac{yK_0}{\tilde{d}^4 h},\,\frac{1}{||\alpha {\tilde{d}}^4 h||}\right\}.
\end{align}
Applying Lemma \ref{TT} and appealing to (\ref{Ogranichenia}), \eqref{Usl3}, and (\ref{W2malkoAsquare}), we get
\begin{equation}\label{W2square_no}
  W_2 \ll y^{\varepsilon}\left(y^{9/10}K +\frac{y^{17/20}K}{q^{1/16}}+
  y^{63/80}K^{15/16}q^{1/16}\right).
\end{equation}
From (\ref{W1intermediate}) and (\ref{W2square_no}), we obtain that for  $y^{3/20}\le D\ll y^{1/5}$, 
\begin{equation}\label{Gama3intermediate}
  W \ll y^{\varepsilon}\left(y^{9/10}K +y^{63/80}K^{15/16}q^{1/16}+\frac{y^{17/20}K}{q^{1/16}}+y^{7/20}K^{1/2}q^{1/2}\right).
\end{equation}

\subsubsection{Estimate of \/ $W$ in the case $1\le D\le y^{3/20}$}\label{new9}

To evaluate the sum $W_1$, we follow the arguments of Section 4.4.1, which led to the estimate \eqref{W_1}. Setting $h=m k \tilde{d}^2$, we have $h\ll MK_0D^2$. From (\ref{Ogranichenia}),(\ref{Usl1}), and Lemma \ref{Trsum2}, we deduce 
\begin{align}\label{W1small}
  W_1  \ll y^{\varepsilon}\sum\limits_{h\ll MK_0D^2}\min\bigg\{\frac{yK_0}{h},\,\frac{1}{||\alpha h||}\bigg\} \ll y^{\varepsilon}\left( \frac{yK}{q}+q+y^{9/10}K  \right).
\end{align}

By partial summation, we obtain the same bound for $W_1'$. 

Reasoning similarly to Section \ref{new8} (see (\ref{W2malkoAsquare}) and (\ref{W21square})), we estimate the sum $W_2$ as follows:
\begin{equation*}\label{W21smallA}
  W_2^2\ll y^{\varepsilon}\Bigg( y^{9/5}K^2 +yDK\sum\limits_{h\ll LD^2K_0}
  \min\bigg\{\frac{yK_0}{h},\,\frac{1}{||\alpha h||}  \bigg\} \Bigg).
\end{equation*}
Applying Lemma \ref{Trsum2} and (\ref{Ogranichenia}), we obtain
\begin{equation}\label{ocW2}
  W_2\ll y^{\varepsilon}\left( y^{9/10}K+\frac{y^{43/40}K}{q^{1/2}}+y^{23/40}K^{1/2}q^{1/2} \right).
\end{equation}
Combining (\ref{W1small}) and (\ref{ocW2}), we conclude that for $1\le D\le y^{3/20}$,
\begin{equation}\label{Gama3small}
  W\ll y^{\varepsilon}\left( y^{9/10}K+\frac{y^{43/40}K}{q^{1/2}}+y^{23/40}K^{1/2}q^{1/2}+q \right).
\end{equation}

\subsubsection{Finalizing the estimate of\/ $\mathcal{U}_3$}
From (\ref{Gama3intermediate}) and (\ref{Gama3small}), we obtain 
  \begin{align}
    \mathcal{U}_3 \ll y^{\varepsilon}  & \bigg(  y^{9/10}K  +\frac{y^{43/40}K}{q^{1/2}} + y^{23/40}K^{1/2}q^{1/2} \nonumber \\
      &{} + \frac{y^{17/20}K}{q^{1/16}}  + y^{63/80}K^{15/16}q^{1/16}+q \bigg). \label{lambda12q_0}
  \end{align}

\subsection{Conclusion of the estimate for \boldmath$\Gamma_2(x)$\unboldmath} \label{conclusion}
From \eqref{new6}, \eqref{veryLargeA0}, \eqref{BigA0}, 
and \eqref{lambda12q_0},
we deduce that the same bound as in \eqref{lambda12q_0} applies to $\Gamma_3(y)$. 
Moreover, from \eqref{Gamma_2-2}, we conclude that the same bound with $y$ replaced by $x$ remains valid for the sum $\Gamma_2(x)$. 

Thus, the estimate of $\Gamma_2(x)$  is now complete. 

\vskip 4pt

\noindent \textbf{Acknowledgements.} The second author was partially supported by the Scientific Fund of Sofia University ``St. Kl. Ohridski'' under Grant 80-10-68/9.4.2024.

\end{document}